\documentclass[12 pt]{article} 
\usepackage[T1]{fontenc} 
\usepackage{ulem} 
\usepackage[utf8]{luainputenc} 
\usepackage{geometry} 
\usepackage[pdftex]{graphicx} 
\usepackage{amstext} 
\usepackage{amssymb} 
\usepackage{amsmath} 
\usepackage{physics}
\usepackage[T1,T2A]{fontenc} 
\usepackage[utf8]{inputenc} 
\usepackage{setspace} 
\usepackage{amsthm}
\usepackage{esint} 
\usepackage{comment} 
\usepackage{babel} 
\usepackage{float} 
\usepackage{fullpage} 
\usepackage{tikz} 
\usepackage{indentfirst} 
\usepackage{mathtext}
\usepackage{ dsfont }
\usepackage{ amsfonts }
\usepackage{ eufrak }
\usepackage{ mathrsfs }
\usepackage{color}
\newtheorem{theorem}{Theorem}
\newtheorem{corollary}{Corollary}
\newtheorem{proposition}{Proposition}

\title{On recurrent properties of Fisher -- Wright's diffusion on $(0,1)$ with mutation}
\author{
R.Yu. Sineokiy\footnote{Moscow State University,  Moscow, Russian Federation; email: ro-ssi @ mail.ru}
\; \& \,
A.Yu. Veretennikov\footnote{Institute for Information Transmission Problems \& National Research University Higher School of Economics, Moscow, Russian Federation; email: ayv @ iitp.ru} 
} 
\begin{document}

\maketitle

\begin{abstract}
One-dimensional Fisher -- Wright diffusion process on the interval $(0,1)$ with mutations is considered. This is a widely known model in population genetics. The goal of the paper is an exponential recurrence of the process, which also implies exponential rate of convergence towards the invariant measure.

~

\noindent
{\bf The keywords}: Wright -- Fisher diffusion; exponential recurrence

~


\end{abstract}

\section{Introduction}
Consider a one-dimensional stochastic differential equation
\begin{equation}\label{e1}
dX_t=[a(1-X_t)-bX_t]dt+\epsilon\sqrt{X_t(1-X_t)}dW_t
\end{equation}
where $a,b,\epsilon>0; X_0=x\in (0,1)$ is a nonrandom initial value of the process, $W_t$ is a one-dimensional Wiener process,  $t \ge 0$.

Such equations were introduced in the populational genetics by S. Wright  and independently by R. Fisher and remain a topical area of investigations until now. For recent sources see, in particular,
\cite{ChenStroock, EpsteinMazzeo} and the references therein.  In particilar, conditions for the existence of the invariant measure and its explicit form in the case without mutations can be found in \cite{ChenStroock}. The paper \cite{DucTranJost18} gives conditions for the existence of the invariant measure (for more general models than (\ref{e1})) and for an exponential rate of convergence; it does not study explicitly recurrence properties, so that the methods are quite different from ours. All issues related to recurrence and convergence rates are interesting not only for populational biology (see \cite{Biol}), but also because this is a sister of the CIR model in stochastic finance, see,  for example, the same source  \cite{DucTranJost18}. For CIR there is a more extensive literature, but we do not use it here, so do not propose further references except for \cite{IlyaGihman} which contains a short proof of sufficiency of Feller's condition for CIR model, which is exploited here for the case of Wright -- Fisher's diffusion.

The equation (\ref{e1}) has a pathwise (and, hence, also weakly) unique strong solution which is pathwise unique according to \cite{YamadaWatanabe}. As it is known -- follows from  uniqueness  \cite{3} -- this solution is a strong Markov process.

It should be highlighted that, although the state space $[0,1]$ is a bounded subset of $\mathbb R^1$, it is still natural to ask about its recurrence properties. Firstly, conditions will be imposed under which the endpoints $0$ and $1$ are not attained; hence, if the initial point $x$ belongs to the open interval, $(0,1)$ then the state space is, actually, not $[0,1]$, but this open interval, which is equivalent to the whole $\mathbb R^1$ topologically. Therefore, it makes sense to talk about recurrence ``from $\pm\infty$'', where the role of $-\infty$ and $+\infty$ are played by $0$ and $1$, respectively. Recurrence is then understood as a return from any small neighbourhood of $0$ or $1$ to a compact set in $(0,1)$. As is well-known, very important features of recurrence are the moments of corresponding hitting times. It will be shown that for the intervals $[a,b]$ with $a$ close to $0$ and $b$ close to $1$, such hitting times have  exponential moments.

The paper is arranged as follows: section 1 is this introduction; section 2 contains main results; section 3 is devoted to the proof of theorem \ref{T1} which is the principal result of this note; section 4 gives a hint about the proof of corollary \ref{C1}, and for the convenience of the readers section 5 provides the proof of proposition \ref{L1}.

\section{Main results}
For any $\alpha \in (0,1/2)$ let
$$
\tau=\tau_\alpha:= \inf(t\ge 0: \; X_t\in [\alpha, 1-\alpha]).
$$

\begin{theorem}\label{T1}
Let Feller's condition
\begin{equation}\label{Fe}
a\wedge b \equiv \min(a,b)> \frac{\epsilon^2}{2}
\end{equation}
be satisfied.
Then for {\bf any} $c>0$ there exist $\alpha>0$, $m>0$ such that for any $x\in (0,1)$

\begin{equation}\label{e3}
\mathbb E_x e^{c\tau_\alpha}\leq C(m)c\alpha^{m+1}((1-x)^{-m}+x^{-m})+1,
\end{equation}
where
$$
C(m)=\max\left(\frac{2}{bm-\epsilon^2 m(m+1)/2}, \frac{2}{am-\epsilon^2 m(m+1)/2}\right) \equiv \frac{2}{(a\wedge b)m-\epsilon^2 m(m+1)/2}.
$$
Moreover,
\begin{equation}\label{e4}
\mathbb E_x \int_{0}^{\tau_\alpha} X_s^{-m-1}ds\leq C(m)c\alpha^{m+1} (x^{-m}+ (1-x)^{-m}).
\end{equation}
\end{theorem}
It will be shown in the proof how to choose $\alpha$ and $m$ given the values of $a,b,c, \epsilon>0$.

\begin{corollary}\label{C1}
Under the assumption (\ref{Fe}) there exists a unique invariant measure $\pi$  for the family of Markov diffusions (\ref{e1}), and there {\bf exist} constants $C,c,m>0$ such that the distribution of $X_t$ -- denoted by $\mu_t^x$ -- converges for any $x\in(0,1)$ to this invariant measure in the total variation metric with the rate admitting the bound \begin{equation}\label{convergence}
\|\mu_t^x - \pi\|_{TV} \le C\exp(-ct)(x^{-m}+ (1-x)^{-m} +2).
\end{equation}
\end{corollary}
Emphasize that unlike in the theorem \ref{T1} here the value $c>0$ may not be arbitrary. The reason for this is that the smaller $\alpha>0$, the greater $c$ can be taken  in (\ref{e3}); however, on the other hand, the smaller $\alpha$, the poorer the nondegeneracy of the diffusion coefficient, which, in turn, decreases the local mixing properties of the process inside the compact $[\alpha,1-\alpha]$. Yet, overall the bound (\ref{e3}) does lead to some exponential estimate for the convergence rate (\ref{convergence}) with respect to time variable, non-uniform with respect to the initial state $x$ from the open interval $(0,1)$.

Recall that the existence and uniqueness of the invariant measure for the equation (\ref{e1}) are well-known; however, both these claims follow again independently and automatically from the theorem \ref{T1} by virtue of the Harris -- Khasminskii principle, see  \cite{7}. Yet, the main news of this corollary is the bound (\ref{convergence}).  

The upper bound (\ref{e3}) leaves open the question of recurrence -- and, hence, of convergence rate -- if the process starts at zero or one, or, if it touches one of these endpoints. (Note that if this happens, it still does not spoil the statement about a strong solution.) So, for the convenience of the reader and having in mind further applications in the next studies let us state the fact of the inattainability of the endpoints of the interval $(0,1)$ under Feller's condition in a proposition; in what follows a brief simple proof will be provided not repeating the original paper \cite{4} (where it is shown rigorously exactly for this model by a different method) but based on
the nice idea from \cite{IlyaGihman} with some slight changes. 

\begin{proposition}\label{L1}
Under the assumptions of theorem \ref{T1} including Feller's condition (\ref{Fe}), 
\begin{equation}\label{eqprop}
\mathbb P_x\left(\exists \, t <\infty: \; X_t(1-X_t) = 0\right) = 0, \quad \forall\, x\in (0,1).
\end{equation}
\end{proposition}

\section{Proof of theorem \ref{T1}}

Let for definiteness $x \in (0,\alpha)$, where $\alpha\in (0,1/2)$ is a constant to be chosen in what follows. The case $x\in (1-\alpha, 1)$ will be commented in the end of this section.  Let $c>0$,  $m>0$; here $c>0$ is arbitrary while  $m>0$ will also be chosen.
Let us consider the function $V(x,t)=\exp(ct)x^{-m}$ and let us compute the stochastic differential of $V(X_t,t)$: by It\^o's formula we have, ($dX_t=[a-(a+b)X_t]dt+\epsilon\sqrt{X_t(1-X_t)}dW_t$)
\begin{align*}
de^{ct}X_t^{-m}=e^{ct}\left\{cdt - m X_t^{-m-1}dX_t + \frac{m(m+1)}2 (dX_t)^2\right\}
 \\\\
=\! e^{ct}\left\{cX_t^{-m}dt \!-\! m X_t^{-m-1} \left([a\!-\!(a\!+\!b)X_t]dt\!+\!\epsilon\sqrt{X_t(1\!-\!X_t)}dW_t\right)
\!+\! \frac{\epsilon^2 m(m\!+\!1)}2 X_t^{-m-1}(1\!-\!X_t)dt\right\}
 \\\\
= e^{ct}X_t^{-m-1}\left\{cX_t^{} - m a + m(a+b)X_t +\! \frac{\epsilon^2 m(m\!+\!1)}2 (1\!-\!X_t)\right\}dt
\!-\! e^{ct}m \epsilon \sqrt{X_t(1-X_t)} X_t^{-m-1}dW_t
\\\\
= \!e^{ct}X_t^{-m-1}\!\left\{\!- m a \!+\! \frac{\epsilon^2 m(m\!+\!1)}2 \!+\! \left(c \!+m(a\!+\!b)\! -\! \frac{\epsilon^2 m(m\!+\!1)}2\right) \!X_t\right\}\!dt\!
- \!e^{ct}m \epsilon \sqrt{X_t(1-X_t)} X_t^{-m-1}dW_t .
\end{align*}
Let us define stopping times
$$
\tau_{t,N}:=\min(\tau,t,\inf(s\ge 0: X_s\wedge (1-X_s)\le 1/N)),
\quad t\ge 0, \; N=3,4,\ldots
$$
Integrating and taking expectations, we obtain,
\begin{align*}
\mathbb E_x e^{c\tau_{t,N}}X_{\tau_{t,N}}^{-m}-x^{-m}
 \\\\
=\mathbb E_x\int\limits_{0}^{\tau_{t,N}} e^{cs}\left[\left(c+(a+b)m-\frac{\epsilon^2}{2}m(m+1)\right)X_s^{-m}+\left(\frac{\epsilon^2}{2}m(m+1)-am\right)X_s^{-m-1}\right]ds.
\end{align*}
Since
$\mathbb E_xe^{c\tau_{t,N}}X_{\tau_{t,N}}^{-m} > 0$, then we may conclude that
\[
x^{-m}\geq \mathbb E_x\int_{0}^{\tau_{t,N}} e^{cs}\left[\left(am-\frac{\epsilon^2}{2}m(m+1)\right)X_s^{-m-1}-(c+(a+b)m-\frac{\epsilon^2}{2}m(m+1))X_s^{-m}\right]ds.
\]
Now let us choose $\alpha>0$ so that
\begin{equation}\label{mult2}
(am-\frac{\epsilon^2}{2}m(m+1))\geq2(c+(a+b)m)\alpha.
\end{equation}
(Here the multiplier $2$ in the right hand side can be replaced by any constant strictly greater than one; any such replacement would change the constants in the resulting inequality but not the principal conclusion.) In turn, for this it suffices that 
\begin{equation}\label{eq:ref1}
\alpha<\frac{(am-\frac{\epsilon^2}{2}m(m+1))}{2(c+(a+b)m)},
\end{equation}
and
\begin{equation}\label{eq:ref3}
am-\frac{\epsilon^2}{2}m(m+1)>0.
\end{equation}
Note that with the multiplier $2$ in (\ref{mult2}) the right hand side in (\ref{eq:ref1}) is automatically less than $1/2$, so that the interval $[\alpha,1-\alpha]$ is not empty. If some other multiplier close to $1$ were used in (\ref{mult2}), it would be natural to take
$$
0< \alpha<\frac{(am-\frac{\epsilon^2}{2}m(m+1))}{2(c+(a+b)m)}
\wedge (\frac12 - \delta)
$$
with any small $\delta \in (0,1/2)$.
The bound (\ref{eq:ref3}) is equivalent to
$$
\frac{2a}{\epsilon^2} > m+1,
$$
and exactly Feller's condition (\ref{Fe}) guarantees that there exists a positive $m$ such that
$$
0<m<\frac{2a}{\epsilon^2} -1.
$$
With this choice of $m>0$ we obtain
\begin{align*}
\mathbb E_x\int_{0}^{\tau_{t,N}} e^{cs}\left[\left(am-\frac{\epsilon^2}{2}m(m+1)\right)X_s^{-m-1}-\left(c+(a+b)m-\frac{\epsilon^2}{2}m(m+1)\right)X_s^{-m}\right]ds
 \\\\
\ge \mathbb E_x\int_{0}^{\tau_{t,N}} e^{cs}\left[\left(am-\frac{\epsilon^2}{2}m(m+1)\right)X_s^{-m-1}-\frac{1}{2}\left(am-\frac{\epsilon^2}{2}m(m+1)\right)X_s^{-m-1}\right]ds
 \\\\
=\frac{1}{2} \mathbb E_x\int_{0}^{\tau_{t,N}} e^{cs}\left(am-\frac{\epsilon^2}{2}m(m+1)\right)X_s^{-m-1}ds.
\end{align*}
Note that the integrand in the last equality is positive due to the condition (\ref{eq:ref3}). Hence, we get
\[
\mathbb E_x\int_{0}^{\tau_{t,N}} e^{cs}\left[\frac{1}{2}\left(am-\frac{\epsilon^2}{2}m(m+1)\right)X_s^{-m-1}\right]ds \leq x^{-m}.
\]
Since $X_s\leq\alpha$, and $e^{cs}\geq1$ for any $s\in[0,\tau_{t,N}]$, $c>0$, then the following two inequalities follow:

\begin{equation}
\mathbb E_x\int_{0}^{\tau_{t,N}}{e^{cs}ds}\leq \frac{2}{am-\epsilon^2m(m+1)/2} \alpha^{m+1} x^{-m}\label{eq:ref4},
\end{equation}
and
\begin{equation}
\mathbb E_x\int_{0}^{\tau_{t,N}} X_s^{-m-1}ds\leq \frac{2}{am-\epsilon^2m(m+1)/2} \, x^{-m}.
\label{eq:ref5}
\end{equation}
Denote
$$
C(m) := \frac{2}{am-\epsilon^2m(m+1)/2}.
$$ Then, by integration of the left hand side in (\ref{eq:ref4}) we obtain for $0< x<\alpha$
\[
\mathbb E_xe^{c\tau_{t,N}}\leq C(m)c\alpha^{m+1} x^{-m}+1.
\]
By virtue of Fatou's lemma (or, by the monotone convergence theorem)  we get as $N\uparrow\infty$ and $t\uparrow\infty$
\[
\mathbb E_xe^{c\tau}\leq C(m)c\alpha^{m+1} x^{-m}+1.
\]
Similarly, applying Fatou's lemma to the inequality (\ref{eq:ref5}), we get as  $N\uparrow\infty$ and $t\uparrow\infty$,

\begin{equation*}
\mathbb E_x\int_{0}^{\tau} X_s^{-m-1}ds\leq C(m) x^{-m}. \label{eq:ref6}
\end{equation*}
This proves the statement of the theorem in the case of $x \in (0, \alpha)$.

The case $x \in (1-\alpha, 1)$ can be treated  similarly via the Lyapunov function $V(t,x) = \exp(ct)(1-x)^{-m}$. Otherwise, the change of variables $y=1-x$ may be used; then the transformation
$$
Y_t := 1-X_t,
$$
reduces the situation $x > 1-\alpha$ to the previous case $x <\alpha$. Here the choice of $m>0$ and $\alpha>0$ has to satisfy the inequalities
\begin{equation*}
bm-\frac{\epsilon^2}{2}m(m+1)>0,
\end{equation*}
and
\begin{equation*}
\alpha<\frac{(bm-\frac{\epsilon^2}{2}m(m+1))}{2(c+(a+b)m)},
\end{equation*}
instead of (\ref{eq:ref1}) and (\ref{eq:ref3}). This means that overall we may take any value of $m$ in the interval
\begin{equation*}
0<m<\frac{2(a\wedge b)}{\epsilon^2} -1,
\end{equation*}
and then any value of $\alpha$ in the interval
\begin{equation*}
0<\alpha<\frac{((a\wedge b)m-\frac{\epsilon^2}{2}m(m+1))}{2(c+(a+b)m)}.
\end{equation*}
The theorem is proved. \hfill QED

\section{Proof of corollary \ref{C1}, sketch}

The statement follows from the inequalities (\ref{e3}) and (\ref{e4}) of the theorem \ref{T1} via one of the standard techniques of the Harris -- Khasminskii principle and coupling, as can be  seen, for example, in  \cite{2,6}. The details will be presented in the next publications for a more general model. \hfill QED

\section{Proof of Proposition \ref{L1}}
NB: This proof is fully independent of the statement of theorem \ref{T1} and of its proof. 

\medskip

\noindent
Let us study the case of the endpoint 0. Denote $B(x):= a - (a+b)x$, and $\sigma(x) := \sqrt{x(1-x)}$. Let 
$$
0< \kappa <\frac{a - \epsilon^2/2}{a+b}, \quad
b_0:= a - (a+b)\kappa, \quad 
n\in \left(0,\frac{2b_0}{\epsilon^2}-1\right]. 
$$
(This value $n$ does not mean to be necessarily integer.)
Note that $b_0>\epsilon^2/2$, so that the choice of such a positive $n$ is possible, and 
$b_0 - (n+1)\epsilon^2/2\ge 0$. 
Consider a Lyapunov function $V(x)=x^{-n}$. By It\^o's formula, 
\[
dX_t^{-n}=[-nB(X_t)X_t^{-n-1}+\frac{n(n+1)\epsilon^2}{2}\sigma(X_t)^2X_t^{-n-2}]dt - n\epsilon \sigma(X_t)X_t^{-n-1}dW_t.
\]
Let any $\beta \in (0, \kappa)$, and define the stopping times
$$
\gamma_\beta:=\min\{t\geq 0: X_t\leq \beta\}, \quad 
\gamma_0:=\min\{t \geq 0: \, X_t=0\},\quad T_{\kappa}=\min\{t\geq 0: X_t \ge \kappa\},
$$ 
and 
$$
\tau_{t,\beta, \kappa}:=\min(\gamma_\beta,T_{\kappa},t) 
\equiv \min(\gamma_\beta \wedge t,T_{\kappa}).
$$ 
It suffices to show that $\mathbb P_x(\gamma_0 < T_{\kappa})=0$. We have, 
\begin{align*}
\mathbb E_x X^{-n}_{\tau_{t,\beta, \kappa}}=x^{-n} + \mathbb E_x\int_{0}^{\tau_{t,\beta, \kappa}} \left[-nB(X_s)X_s^{-n-1}+\frac{n(n+1)\epsilon^2}{2}\sigma(X_s)^2X_s^{-n-2}\right]\,ds
 \\\\
\le x^{-n} + \mathbb E_x\int_{0}^{\tau_{t,\beta, \kappa}}\left[-n b_0 X_s^{-n-1}+\frac{n(n+1)\epsilon^2}{2}X_s(1-X_s)X_s^{-n-2}\right]ds
 \\\\
= x^{-n} + \mathbb E_x\int_{0}^{\tau_{t,\beta, \kappa}}\left[-nX_s^{-n-1} \left(b_0 -\frac{(n+1)\epsilon^2}{2}\right) - \frac{n(n+1)\epsilon^2}{2}X_s^{-n}\right]ds
 \\\\
\le x^{-n} - \mathbb E_x\int_{0}^{\tau_{t,\beta, \kappa}}{\frac{n(n+1)\epsilon^2}{2}X_s^{-n}}ds\leq x^{-n},
\end{align*}
the one before the last inequality because $\displaystyle b_0-\frac{(n+1)\epsilon^2}{2}\ge 0$ by the choice of $n$. Now by the Bienaym\'e -- Chebyshev -- Markov inequality 
\[ 
\mathbb P_x(X_{\tau_{t,\beta, \kappa}}\le \beta)=\mathbb P_x(X^{-n}_{\tau_{t,\beta, \kappa}}\ge \beta^{-n})\leq \beta^n x^{-n}\xrightarrow{} 0, \quad \beta \rightarrow 0.
\]
So for any $t\ge 0$, 
\begin{align*}
\mathbb P_x(\gamma_\beta\wedge t < T_{\kappa}, \gamma_\beta  \le t) = \mathbb P_x(X_{\tau_{t,\beta,\kappa}} \le \beta, \gamma_\beta\wedge t < T_{\kappa}, \gamma_\beta \le t) 
\le \mathbb P_x(X_{\tau_{t,\beta,\kappa}}\le \beta) 
\to 0, \quad \beta \to 0.
\end{align*}
Since $\gamma_0>\gamma_\beta$ for any $\beta$, we obtain for each $t$:
\[
\mathbb P_x(\gamma_0 \wedge t < T_{\kappa}, \gamma_0 \le t)\le  \mathbb P_x(\gamma_\beta\wedge t < T_{\kappa}, \gamma_\beta  \le t) 
\to 0, \quad \beta \rightarrow 0. 
\]
Hence, $\mathbb P_x(\gamma_0 < T_{\kappa})=0$, and so also 
$$
\mathbb P_x\left(\exists \, t <\infty: \; X_t = 0\right) = 0, \quad \forall\, x\in (0,1), 
$$
as required.
The case of the endpoint  1 is considered similarly; these two cases combined finally lead to the equation (\ref{eqprop}). \hfill QED

\section*{Acknowledgements}

For the second author this study in part of proposition \ref{L1}  was prepared within the framework of the HSE University Basic Research Program, and in part of corollary \ref{C1} it was funded by the Russian Science Foundation grant  17-11-0198 (extended).


\begin{thebibliography}{99}

\bibitem{ChenStroock}
L. Chen and D. W. Stroock,
The fundamental solution to the
Wright--Fisher equation,
{\it SIAM J. Math. Anal.},
{\bf42}(2010), no. 2, , 539--567.
https://doi.org/10.1137/090764207

\bibitem{DucTranJost18}
L. H. Duc, T. D. Trana and J. Jost,
Ergodicity of scalar stochastic differential equations with Hölder continuous coefficients,
{\it Stochastic Processes Appl.},
{\bf128}(2018), no. 10, , 3253--3272.
https://doi.org/10.1016/j.spa.2017.10.014

\bibitem{EpsteinMazzeo}
C. L. Epstein and R. Mazzeo,
Wright–Fisher diffusion in one dimension,
{\it SIAM J. Math. Anal.}, {\bf42}(2010), no. 2,  568-–608.
https://doi.org/10.1137/090766152

\bibitem{4}
W. Feller, Two singular diffusion problems, {\it Ann. of Math.}, {\bf54}(1951), 173--182. https://doi.org/10.2307/1969318

\bibitem{IlyaGihman}
Ilya I. Gikhman, A short remark on Feller’s square root condition,  https://ssrn.com/abstract=1756450 (2011); https://dx.doi.org/10.2139/ssrn.1756450

\bibitem{7}
R. Z. Khas’minskii,
Ergodic properties of recurrent diffusion processes and stabilization of the solution to the Cauchy problem for parabolic equations,
{\it Theory Probab. Appl.,} {\bf5}(1960), no. 2 , 179--196.
https://doi.org/10.1137/1105016

\bibitem{3}
N. V. Krylov,
On the selection of a Markov process from a system of processes and the construction of quasi-diffusion processes,
{\it Math. USSR-Izv.,} {\bf7}(1973), no. 3, 691--709.
http://dx.doi.org/10.1070/IM1973v007n03ABEH001971

\bibitem{Biol}
M. Steinr\"ucken, R. Wang and  Y. S. Song,
An explicit transition density expansion for a multi-allelic Wright-Fisher diffusion with general diploid selection,
{\it Theor. Popul. Biol.}, {\bf83}(2013),  1--14. doi: 10.1016/j.tpb.2012.10.006.

\bibitem{2}
A. Yu. Veretennikov, On polynomial mixing bounds for stochastic differential equations, {\it Stochastic Processes Appl.,} {\bf70}(1997), 115--127. https://doi.org/10.1016/S0304-4149(97)00056-2

\bibitem{6}
A. Yu. Veretennikov,
On polynomial mixing and convergence rate for stochastic difference and differential equations,
{\it Theory Probab. Appl.}, {\bf44}(2000), no. 2, 361–-374.
https://doi.org/10.1137/S0040585X97977550

\bibitem{YamadaWatanabe}
T. Yamada and S. Watanabe, On the uniqueness of solutions of stochastic differential equations, {\it J. Math. Kyoto Univ.}, {\bf11}(1971), no. 1, 155--167. DOI: 10.1215/kjm/1250523691

\end{thebibliography}
\end{document}